\documentclass{amsart}
\usepackage{amsfonts}
\usepackage{amsthm}
\usepackage{amsmath}
\usepackage{amsfonts}
\usepackage{latexsym}
\usepackage{amssymb}
\usepackage{verbatim}

\newcommand{\U}{\mathcal{U}}

\newcommand{\RR}{\mathbb{R}}
\newcommand{\NN}{\mathbb{N}}

\newcommand{\FF}{\mathbb{F}}

\newtheorem{fed}{Definition}[section]
\newtheorem{teor}[fed]{Theorem}
\newtheorem*{teo*}{Theorem}
\newtheorem{lema}[fed]{Lemma}
\newtheorem{coro}[fed]{Corollary}

\theoremstyle{definition}
\newtheorem{rema}[fed]{Remark}

\newtheorem{exa}[fed]{Example}
\newtheorem{exas}[fed]{Examples}

\DeclareMathOperator{\Tr}{Tr} \DeclareMathOperator{\tr}{tr}

\newcommand{\CC}{\mathbb{C}}

\begin{document}

\title[  Optimal systems for erasures and the
q-potential ] {Optimal reconstruction systems for erasures and for
the q-potential}

\author[Pedro G. Massey]{Pedro G. Massey* }

\thanks{*Depto. Matem\'atica-FCE-UNLP and Instituto Argentino de Matemática-CONICET, Argentina }
\thanks{ Partially supported by CONICET , Universidad Nacional de La Plata and ANPCYT}

\begin{abstract}
We introduce the $q$-potential as an extension of the
Benedetto-Fickus frame potential, defined on general reconstruction
systems and we show that protocols are the minimizers of this potential
under certain restrictions. We extend recent results of B.G. Bodmann
on the structure of optimal protocols with respect to 1 and 2 lost
packets where the worst (normalized) reconstruction error is
computed with respect to a compatible unitarily invariant norm. We
finally describe necessary and sufficient (spectral) conditions,
that we call $q$-fundamental inequalities, for the existence of
protocols with prescribed properties by relating this problem to
Klyachko's and Fulton's theory on sums of hermitian operators.
\end{abstract}

\maketitle

\noindent{\bf Keywords. Reconstruction systems; $q$-potential;
Erasures; Compatible unitarily invariant norm; $q$-fundamental
inequality.}


\section{Introduction}

Signal transmission through a noisy channel, such as digital media
through the Internet, typically uses the following strategy: a
generic signal is decomposed (encoded) into a sequence of
coefficients which are then grouped into a number of packets of the
same size. These packets are then sent through the noisy channel.
For practical purposes, we shall assume that the noise in the
channel cannot affect the integrity of the data in each packet; we
can think that these \emph{small} pieces of data are protected by an
efficient error-correcting algorithm. Still, the noise of the
channel may cause such a delay in the transmission or even the loss
of some packets that the reconstruction of the signal is done
possibly without the whole set of packets. At this point there are
at least two different procedures to follow: we can attempt to
reconstruct the signal by altering the decoding algorithm (taking
into account the fact that there are some lost packets) or we can
apply the decoding process to the received packets and accept that
we obtain only an approximation of the original signal.

 In the present paper we will adopt the second alternative
 for
 the reconstruction of the signal. Hence we search for
 encoding-decoding schemes that minimize, with respect to some measure,
 the worst case error between (a normalization of) the original  signal and
 the reconstructed signal for a fixed number of packet losses, under certain
 hierarchies (see the beginning of Section \ref{sec 4} for a description of these hierarchies).  This and similar problems have been considered recently by Casazza and Kovacevic \cite{Cas},  Heath and Stromer \cite{Stro},
Holmes and Paulsen \cite{HolPau},
Bodmann and Paulsen \cite{BodPau} and Bodmann \cite{Bod} (and on the
related work of Bodmann, Paulsen and Kribs \cite{Pau}) where they
describe the structure of optimal encoding-decoding schemes based on
a particular choice to measure the worst case reconstruction error.
Some of the results in the present paper can then be described as
generalizations of some of the results obtained in those works, as
we show that the previously mentioned optimal schemes are actually
optimal with respect to a continuous family of measures (that
includes most of the typical choices) of the worst case
reconstruction error in the more general setting of
block-(encoding-decoding) introduced in \cite{Bod}.
 Our approach and techniques related with these problems are derived
here as a generalization of those in \cite{Bod}.

The  optimal schemes found in the frame-based transmission model
(under suitable restrictions) are related with the so-called
Parseval (or more generally tight) frames. These frames are also
important for applications since they allow for representations of
signals that are formally the same as those given by an orthonormal
basis, but with the additional property of redundancy. There is a
natural generalization of Parseval frames introduced by Bodmann in
\cite{Bod}, the so-called \emph{protocols}, which is the starting
point for the development of the theory of optimal protocols under
packet-erasures in that paper. In this setting, the optimal
protocols correspond to some \emph{projective} protocols, which were
originally introduced by Casazza and Kutyniok \cite{Caspalo} under
the name of (Parseval) frames of subspaces, and recently have also
been called fusion frames \cite{Cas2}. But there are more general
\emph{reconstruction systems} (see Definition \ref{defi recons}) than
protocols, just as there are more general frames than Parseval
frames.

In order to investigate possible advantages of general protocols in
the class of reconstruction systems we introduce what we called the
$q$-potential, which is a generalization of the frame potential
defined by Benedetto and Fickus  in \cite{BF} and further considered
in \cite{FJKO} and \cite{casazza2}. In our case the $q$-potential of
a reconstruction system takes values in the cone of positive
matrices, rather than numerical values, a fact that makes it difficult
to compare $q$-potentials of different systems. Still, we show that
under suitable conditions, protocols are the minimizers of the
$q$-potential within reconstruction systems with
respect to \emph{(sub)majorization} and thus we obtain lower bounds
and minimizers of a family of (anti)entropic measures of the
$q$-potential. These results indicate that protocols are indeed a
good stating point for the theory of block-erasures.

Although there are interesting techniques to
construct 2-uniform protocols i.e. protocols that are optimal for
two packet losses (see \cite{Bod}, \cite{BodPau}, \cite{HolPau}),
the problem of finding necessary and sufficient conditions for the
existence of protocols that are optimal for one packet loss has been
considered open (see the discussion in \cite{Pau}). We relate this
problem to a problem solved by Klyachko \cite{Klya} and Fulton
\cite{Ful} related with Horn's conjecture on the sums of hermitian
matrices and hence we obtain a characterization of the existence of
such optimal protocols. This result can be regarded as an extension
of the equivalence of the Schur-Horn problem on the main diagonal of
an hermitian operator with prescribed spectrum and the problem of
finding necessary and sufficient conditions for the existence of a
frame for a finite dimensional Hilbert space with prescribed norms
and frame operator as described in \cite{JDM} (see also
\cite{casazza3,{Lar},{Tropp}}), using the notion of
extended majorization as described in \cite{Mas}. We then derive
derive the $q$-fundamental inequality (see Corollary \ref{q fund ineq}), that
is a generalization of the fundamental inequality found in
\cite{casazza2}.

The  paper is organized as follows. In section \ref{spreliminares},
after introducing some notation, we recall the notion of
majorization and some of its properties. Then we distinguish a class
of unitarily invariant norms (that we call compatible) that plays a
key role here. We end the preliminaries by briefly describing
 the basic elements of Klyachko's theory on the sums of
hermitian operators. We begin section \ref{sec 3} by introducing the
$q$-potential defined in the class of reconstruction systems and
show that the protocols are the minimizers of this positive operator
function with respect to submajorization. Thus, it is natural to
restrict the analysis of optimal reconstruction systems for erasures
to protocols. In subsection \ref{subsec 4.1} we give a complete
description of optimal protocols for one packet loss, when we base
the measure of the worst case reconstruction error on a compatible
unitarily invariant norm. In subsection \ref{subsec 4.2} we deal
with the case of two lost packets where we show explicitly a family
of optimal protocols, when restricted to a certain family of optimal
protocols for one loss packet. We then show that this restriction is
automatically satisfied by optimal frames for one coefficient loss
and obtain a generalization of previous results on the structure of
optimal frames for two lost packets. Finally, in section \ref{sec 5}
we consider the problem of designing protocols with prescribed
additional properties. In particular, we find necessary and
sufficient conditions for the existence of optimal protocols for one
packet loss, in terms of a finite set of linear inequalities.

\medskip

\noindent  {\bf Acknowledgments.} I would like to thank the members
of the `` Mischa Cotlar '' seminar at the Instituto Argentino de
Matem\'atica - CONICET, who listened to an early version of this
work and made useful suggestions regarding the material herein.

\section{Preliminaries}\label{spreliminares}

In this note we shall denote by $\mathcal H=\mathbb F^d$ and
$\mathcal K=\mathbb F^l$, where $\mathbb F$ stands for $\RR$ or
$\CC$ and $l\leq d$. Hence, if $l<d$ there is a natural injection
$\iota:\mathcal K\rightarrow \mathcal H$ such that
$\iota(x)=(x,0_{d-l})$, where $0_{d-l}$ denotes the zero vector in
$\RR^{d-l}$. Moreover, $\mathcal H=\mathcal K\oplus \mathcal K'$
under the identification given by $\iota$, with $\mathcal K'=\iota(K)^\perp$. Given
$r,t\in \NN$ we denote by $M_{r,t}(\FF)$ the $\FF$-algebra of all
$r\times t$ matrices with entries in $\FF$. For simplicity we note
$M_r(\FF)$ instead of $M_{r,r}(\FF)$. We further consider
$M_r(\FF)^{sa}$, $M_{r}(\FF)^+$ and $\U(r)$ that are the real space
of self-adjoint matrices, the cone of positive semi-definite
matrices and the group of unitary matrices over $\FF$, respectively.
If $A\in M_d(\CC)^{sa}$ then we denote by $\lambda(A)\in \RR^d$ the
vector of eigenvalues of $A$ (counting multiplicities) with its
entries arranged in decreasing order.
 By fixing the canonical basis in $\mathcal
H$ and $\mathcal K$ respectively, we shall identify $\mathcal
L(\mathcal H)$, $\mathcal L(\mathcal K)$ and $\mathcal L(\mathcal
H,\mathcal K)$ with $M_{d}(\FF)$, $M_{l}(\FF)$ and $M_{l,d}(\FF)$
respectively. The vector $\mathbf e_d\in \RR^d$ is the vector with
all its entries equal to 1.  Finally, if $X$ is a finite set then
$|X|$ denotes its cardinal.

\subsection{Submajorization in $M_l(\CC)^{sa}$}

Given $x\in \RR^l$ we denote by $x^\downarrow\in \RR^l$ the vector
obtained by re-arrangement of the coordinates of $x$ in
non-increasing order. Given $x,\,y\in \RR^l$ we say that $x$ is
\emph{submajorized} by $y$, and write $x\prec_w y$ if
\begin{equation}\sum_{i=1}^k x^\downarrow _i\leq \sum_{i=1}^k y^\downarrow _i \, , \ \ \text{
for } \ 1\leq k\leq l.
\end{equation} If we further have that
$\tr(x):=\sum_{i=1}^lx_i=\sum_{i=1}^l y_i$ then we say that $x$ is
majorized by $y$, and write $x\prec y$.

\begin{exa}\label{ejem}
As an elementary example, that we shall use repeatedly, let $x\in \RR_{\geq 0}^l$ and $0\leq a\leq \tr(x)\leq b$:
then the reader can easily verify that
\begin{equation}\label{eq exa}\frac{a}{l}\,\mathbf e_l\prec_w x\prec_w
b\, e_1
\end{equation}
\end{exa}
The following result, that we shall need in the sequel, is an slight
strengthening of the previous example.
\begin{lema}\label{lemmayo}
Let $\alpha^\downarrow=(\alpha_1,\alpha_2),\,\beta^\downarrow=(b_1\,
\mathbf e_l, b_2\,\mathbf e_l)\in \RR_{\geq 0}^{2l}$ be such that
$\tr(\alpha)\geq \tr(\beta)$ and $\tr(\alpha_1)\geq b_1\,l$. Then
$\beta$ is submajorized by $\alpha$.
\end{lema}
\begin{proof}
Since $\tr(\alpha_1)\geq b_1\,l$ then by Example \ref{ejem}, $b_1\,\mathbf e_l\prec_w
\alpha_1=(a_1^{(1)},\ldots,a_l^{(1)})$. Hence, if $1\leq k\leq l$
then $\sum_{i=1}^k\alpha^\downarrow_i=\sum_{i=1}^k a_i^{(1)}\geq k\,
b_1=\sum_{i=1}^k \beta^\downarrow_i$. If
$\alpha_2=(a_1^{(2)},\ldots, a_l^{(2)})$ then define
$\gamma=(a_1^{(2)}+ (\tr(\alpha_1)-b_1\,l),a_2^{(2)},\ldots, a_l^{(2)})$ and
note that $\gamma=\gamma^\downarrow\in \RR_{\geq 0}^l$. Since
$\tr(\gamma)=\tr(\alpha_1)+\tr(\alpha_2)-b_1\,l\geq b_2\,l$ then we
conclude again that $b_2\,\mathbf e_l\prec_w \gamma$. If $1\leq k\leq l$
then $\sum_{i=1}^k\gamma_i=\sum_{i=1}^{l+k}\alpha^\downarrow_i
-b_1\,l\geq b_2\,k$ and the lemma follows from this last fact.
\end{proof}

(Sub)majorization between vectors is extended by T. Ando in \cite{Ando} to (sub)majorization
between self-adjoint matrices as follows : given $A,\,B\in
M_l(\CC)^{sa}$ then we say that $A$ is submajorized by $B$, and
write $A\prec _w B$, if $\lambda(A)\prec_w \lambda(B)$. If we
further have that $\tr(A)=\tr(B)$ then we say that $A$ is majorized
by $B$ and write $A\prec B$.

Although simple, submajorization plays a central role in
optimization problems with respect to convex functionals and
unitarily invariant norms, as the following result shows (for a
detailed account in majorization see Bhatia's book \cite{Ba}).

\begin{teor}\label{props submayo}
Let $A,\,B\in M_l(\FF)^{sa}$. Then, the following statements are
equivalent:
\begin{enumerate}
\item $A\prec_w B$.
\item For every unitarily invariant norm $\|\cdot\|$ in $M_l(\FF)$ we have
$\|A\|\leq \|B\|$.
\item For every increasing convex function $f:\RR\rightarrow \RR$ we
have $\tr(f(A))\leq \tr(f(B))$.
\end{enumerate}Moreover, if $A\prec_w B$ and there exists an
increasing  strictly convex function $f:\RR\rightarrow \RR$ such
that $\tr(f(A))=\tr(f(B))$ then there exists $U\in \U(l)$ such that
$A=U^*BU$.
\end{teor}

Recall that given a unitarily invariant norm (henceforth abbreviated u.i.n.) $\|\cdot\|$ in
$M_l(\CC)$  there exists an
associated symmetric gauge function $\psi:\RR^l\rightarrow \RR_{\geq 0}$ such
that $\|A\|=\psi(s(A))$, where $s(A)=\lambda(|A|)\in \RR^l$ is the
vector of singular values of $A$. Next we describe a particular
class of u.i.n's that we shall consider in the sequel.
\begin{fed}\label{defi de uin comp}
A sequence $\{ \|\cdot\|_n\}_n$ such that for each $n\in \NN$
$\|\cdot\|_n$ is a u.i.n. in $M_n(\FF)$ is compatible if 
\begin{equation}\label{eccompnorm}\left \|\begin{pmatrix}X & 0 \\ 0 & \mathbf 0_{t}
\end{pmatrix}\right\|_{r+t}=\|X\|_r\end{equation} for every
$X\in M_r(\FF)$, where $\mathbf 0_t\in M_t(\RR)$ is the
zero matrix. If $\psi_n$ is the symmetric gauge function associated
with $\|\cdot\|_n$ then \eqref{eccompnorm} is equivalent to
$\psi_{r+t}(x,0_{t})=\psi_{r}(x)$, where  $x\in \RR^r$ and $0_t\in
\RR^t$ is the zero vector. In this case, we simply write $\|\cdot\|$
and $\psi$ respectively to denote the norms and functions of any
order.
\end{fed}

Let $V:\mathcal H\rightarrow \mathcal K$ be a linear operator and assume that dim$\mathcal H=d>l=$dim$\mathcal  K$. Then,
it is well known that there exists a unitary operator $U\in
\U(d)$ such that
\begin{equation*}\label{relac vv}
U^* \begin{pmatrix} VV^* & 0 \\ 0 & \mathbf 0_{d-l}\end{pmatrix}\,U=V^*V
\end{equation*}where the above block matrix representation is with
respect to the decomposition $\mathcal H=\mathcal K\oplus \mathcal
K'$ as described in the preliminaries. Hence, if $\|\cdot\|$ is a compatible u.i.n. in the sense of
definition \ref{defi de uin comp} we have that $\|VV^*\|=\|V^*V\|$.
This last equality is our main motivation to consider these norms.

We shall use systematically the following facts, which are an
elementary consequence of the previous results: if $\|\cdot\|$ is an
arbitrary u.i.n. in $ M_l(\FF) $ with associated symmetric gauge
function $\psi$ then for every $A\in M_l(\FF)^{+} $ (resp. $x\in
\RR_{\geq 0}^l$) we have \begin{equation*}\|A\|\geq
\frac{\tr(A)}{l}\,\|I_l\|=\tr(A)\,\eta_\psi(l) \ \ \ (\text{resp.} \
\psi(x)\geq \frac{\tr(x)}{l} \psi(\mathbf
e_l)=\tr(x)\,\eta_\psi(l))\end{equation*} where $\eta_\psi(l)=\frac{\|I_l\|}{l}=\frac{\psi(\mathbf e_l)}{l}$, since $\frac{\tr(A)}{l}\,\mathbf e_l \prec\lambda(A)$ and $\frac{\tr(x)}{l}\,\mathbf e_l\prec x$ respectively.

\begin{fed}\label{defi de strict}
A compatible u.i.n. $\|\cdot\|$ is strict if, for any $A\in
M_l(\FF)^+$ then
$$\|A\|=\tr(A)\, \eta_\psi(l) \ \Rightarrow \ A=\frac{\tr(A)}{l} \
I,$$ where $\psi$ is the symmetric gauge function associated with
$\|\cdot\|$ and $\eta_\psi(l)=\frac{\psi(\mathbf e_l)}{l}$.
Equivalently, $\|\cdot\|$ is strict if for $x\in \RR_{\geq 0}^l$
such that $\psi(x)=\tr(x) \,\eta_\psi(l)$ then $x=\frac{\tr(x)}{l}\
\mathbf e_l$.
\end{fed}

\begin{exas}\label{comments before}
As an example of compatible unitarily invariant norm, let us
consider the $p$-norms $\|\cdot \|_p$, with $1\leq p\leq \infty$. On
the other hand, if $1< p\leq \infty$ then $\|\cdot\|_p$ is an strict
norm. Moreover, if $1<p<\infty$ then $f_p(x)=x^p$ is an strictly
convex function and hence the following stronger property holds (see
Theorem \ref{props submayo}): if $A,\,B\in M_l(\CC)^{sa}$ are such
that $A\prec_w B$ and $\|A\|_p=\|B\|_p$ then $A=U^*BU$ for some
$U\in \U(l)$.
\end{exas}

\subsection{Klyachko's and Fulton's spectral theory on sums of hermitian matrices}

In what follows we describe the basic facts about the spectral
characterization of the sums of hermitian matrices obtained by
Klyachko \cite{Klya} and Fulton \cite{Ful}.

Let $\mathcal S_r^d=\{(j_1,\ldots,j_r):\ 1\leq j_1<j_2\ldots<j_r\leq
d\}$. For $J=(j_1,\ldots,j_r)\in \mathcal S_r^d$, define the
associated partition
$$\lambda(J)=(j_r-r,\ldots,j_1-1).$$Denote by $LR_r^{\,d}(m)$ the set of
$(m+1)$-tuples $(J_0,\ldots,J_m)\in (\mathcal S_r^d)^{(m+1)}$, such
that the Littlewood-Richardson coefficient of the associated
partitions $\lambda(J_0),\ldots,\lambda(J_m)$ is positive, i.e. one
can generate the Young diagram of $\lambda(J_0)$ from those of
$\lambda(J_1),\ldots,$ $\lambda(J_m)$ according to the
Littlewood-Richardson rule (see \cite{Ful0}). With these notations
and terminologies we have
\begin{teor}\label{teoK}Let
$\lambda_i=\lambda_i^\downarrow=(\lambda^{(i)}_1,\ldots,\lambda^{(i)}_d)\in
\RR^d$ for $i=0,\ldots,m$. Then, the following statements are
equivalent:
\begin{enumerate}
\item There exists $A_i\in M_d(\CC)^{sa}$ with
$\lambda(A_i)=\lambda_i$ for $0\leq i\leq m$ and such that
$$A_0=A_1+\ldots+A_m.$$
\item For each $r\in \{1,\ldots,d\}$ and $(J_0,\ldots,J_m)\in
LR_r^{\,d}(m)$ we have \begin{equation}\label{comp ec} \sum_{j\in
J_0}\lambda^{(0)}_j\leq\sum_{i=1}^m\sum_{j\in J_i}\lambda^{(i)}_j
\end{equation} plus the condition $\sum_{j=1}^d\lambda^{(0)}_j=
\sum_{i=1}^m\sum_{j=1}^d\lambda^{(i)}_j$.
\end{enumerate}
\end{teor}
We shall refer to the inequalities in \eqref{comp ec} as
\emph{Klyachko's compatibility inequalities}.

For comments on further developments related with the previous
theorem see Remark \ref{comentarios LiPoon}

\section{ Optimality of  $\mathbf{ (m,l,d)}$-protocols for the
q-potential}\label{sec 3}

In what follows we consider $(m,l,d)$-reconstruction systems, which
are more general system of operators than those considered in
\cite{BF}, \cite{Bod}, \cite{Pau}, \cite{BodPau},  \cite{HolPau} and
\cite{MR}, that also have an associated reconstruction algorithm. In
what follows $\mathcal H$ and $\mathcal K$ denote (real or complex)
Hilbert spaces of dimensions $d$ and $l$ respectively, with $l< d$.

\begin{fed}\label{defi recons}
A family $\{V_i\}_{i=1}^m$  is an $(m,l,d)$-reconstruction system if
for $1\leq i\leq m$ $V_i:\mathcal H\rightarrow \mathcal K$ and are
such that $\sum_{i=1}^m V_i^*V_i=S$ is an invertible (positive)
operator.
\end{fed}

 Notice that an $(m,1,d)$-reconstruction system is a frame
 \cite{casart} in the usual sense.

Recall that an $(m,l,d)$-protocol on the Hilbert space $\mathcal H$
\cite{Bod} is a family $\{V_i\}_{i=1}^m$ such that $V_i:\mathcal
H\rightarrow \mathcal K$ for $1\leq i\leq m$ and
$\sum_{i=1}^mV_i^*V_i=I_d$ (see also \cite{Pau}, where protocols are
related to $C^*$-encodings with noiseless subsystems). Clearly,
$(m,l,d)$-protocols are $(m,l,d)$-reconstruction systems in the
sense of definition \ref{defi recons}.

If $\{V_i\}_{i=1}^m$ is an $(m,l,d)$-reconstruction system  then we
consider its \emph{analysis operator} $V:\mathcal H\rightarrow
\oplus_{i=1}^m \mathcal K$ given by $Vx=\oplus_{i=1}^m V_ix$ ;
similarly, we consider its \emph{synthesis operators} given by $V^*$
i.e. $V^*\oplus_{i=1}^m y_i=\sum_{i=1}^m V_i^*y_i$.  For a general
$(m,l,d)$-reconstruction system $\{V_i\}_i$ such that $\sum_{i=1}^m
V_iV_i^*=S$ we have
\begin{equation*}\label{hecho1}
\sum_{i=1}^m V_i^*S^{-1}V_i=I_d
\end{equation*} and thus, we obtain the reconstruction formula
\begin{equation*}\label{ec recons}
x= \sum_{i=1}^m V_i^*S^{-1}(V_i \,x).
\end{equation*} In this context $S=V^*V$ is called the
 \emph{reconstruction system operator} of $\{V_i\}_{i}$ while $G=VV^*$ is called the
 \emph{Grammian operator} of for $\{V_i\}_{i}$.
It is easy to see that in this case $\{S^{-1}V_i\}$ is also an
$(m,l,d)$-reconstruction system, that we call  the \emph{ dual}
reconstruction system associated to $\{V_i\}_i$. For practical
purposes, an encoding-decoding scheme based on the
$(m,l,d)$-reconstruction system above involves the problem of
inverting the reconstruction system operator $S$. One of the
advantages of considering $(m,l,d)$-protocols for applications is
that the reconstruction system operator in this case is $I_d$. As we
shall see $(m,l,d)$-protocols are optimal in other senses, too.

\medskip

In the seminal work \cite{BF} Benedetto and Fickus introduced the
so-called frame potential, as a potential function for the frame
force. The structure of minimizers of the frame potential under
several restrictions \cite {BF}, \cite{FJKO}, \cite{casazza2} and
\cite{MR} have been obtained, since these are considered as stable
configurations with respect to the frame force. This has motivated possible physical
interpretations of families of frames, such as (uniform) tight frames
\cite{casazza2}. Moreover, in \cite{MR} it is shown that the
minimizers of the frame potential (under suitable restrictions) have
structural properties implying their stability with respect to a
more general family of convex functionals that contains the frame
potential of Benedetto and Fickus.

In what follows we introduce the $q$-potential of a reconstruction
system (regardless of an underlying force inducing this potential),
which is a positive semi-definite matrix. Then, we consider two
optimization problems associated with this potential (see Theorems
\ref{propopot} and \ref{infi} below).

\begin{fed} Let $\{V_i\}_{i=1}^m$ be an $(m,l,d)$-reconstruction system on the Hilbert space $\mathcal H$. Then, the $q$-potential of the reconstruction system is defined as \begin{equation*}\label{defi del qpot}
\mathbf P_q(V)=\sum_{i,j=1}^m |V_iV_j^*|^2 \in M_l(\CC)^+.
\end{equation*}
\end{fed}It is straightforward that
the $q$-potential above is the value $\Tr_m((VV^*)^2)\in M_l(\CC)$
i.e. the partial trace of the square of the Grammian operator $VV^*$
with respect to the decomposition $M_{m\cdot
l}(\FF)=M_{m}(M_l(\FF))$. Note that the $q$-potential coincides with
the Benedetto-Fickus potential in the case $l=1$. In contrast to the
Benedetto-Fickus potential, there is no natural way a priori to
compare the $q$-potential of two $(m,l,d)$-reconstruction systems
when $l>1$.

In order to state the following result we recall some distinguished classes of protocols. We say that an $(m,l,d)$-protocol $\{V_i\}_i$ is {\it projective} if for each $1\leq i\leq m$ then $V_i^*V_i=w_i\,P_i$, where $P_i$ is an orthogonal projection in $M_d(\CC)$ and $w_i>0$ are called the associated weights. If the weights of a projective $(m,l,d)$-protocol are equal then we say that it is uniformly weighted (and we abbreviate this by u.w.p). Finally, we say that an $(m,l,d)$-protocol is {\it rank-$l$}, if rank$(V_i^*)=l$ for $1\leq i\leq m$. Notice that if $\{V_i\}_i$ is a rank-$l$ projective $(m,l,d)$-protocol then $V_iV_i^*=w_i\,I_l$ with $w_i>0$, for $1\leq i\leq m$.

\begin{teor}[Optimality of general protocols]\label{propopot}
Let $\{V_i\}_{i=1}^m$ be an $(m,l,d)$-reconstruction system on the
Hilbert space $\mathcal H$ such that $ \tr(V^*V)=\sum_{i=1}^m
\tr(V_i^*V_i)\geq d$. Then, \begin{equation}\label{hora de comer}
\frac{d}{l}\,I_l\prec_w \mathbf P_q(V)
\end{equation}Hence,
for every u.i.n. $\|\cdot\|$ on $M_l(\CC)$ with associated symmetric gauge
function $\psi$ we have
\begin{equation}\label{desiuin} d\cdot \eta_\psi(l)\leq \|\mathbf
P_q(V)\|
\end{equation}
and for every increasing convex function $f:\RR_{\geq 0}\rightarrow
\RR$ with $f(0)=0$ we have \begin{equation}\label{ultimomomento}
l\cdot f(\frac{d}{l})\leq \tr(f(\mathbf P_q(V))).
\end{equation}
If majorization holds in \eqref{hora de comer} or there exists
u.i.n. $\|\cdot\|$ such that equality holds in \eqref{desiuin} or if
there exists an increasing  strictly convex function $f:\RR_{\geq
0}\rightarrow \RR_{\geq 0}$ with $f(0)=0$ such that equality holds in
\eqref{ultimomomento} then $\{V_i\}$ is an $(m,l,d)$-protocol.

Moreover,  if $\{V_i\}_i$ is a projective rank-$l$ $(m,l,d)$-protocol
then majorization holds in \eqref{hora de comer} and the lower
bounds in \eqref{desiuin} and \eqref{ultimomomento} are attained for
each u.i.n. and each function as above, respectively.
\end{teor}

\begin{proof}
Since $\tr(V^*V)\geq d$ then it follows that $I_d\prec_w V^*V\in
M_d(\CC)$ and thus $d=\tr(I_d ^2)\leq \tr((V^*V)^2)=\tr((VV^*)^2)$.
Hence,
\begin{equation}\label{tempranillo} d\leq \tr((VV^*)^2)=\tr(\mathbf
P_q(V))\ \ \Rightarrow \ \ \frac{d}{l}I_l\prec_w \mathbf P_q(V)\in
M_l(\CC).
\end{equation}Notice that \eqref{desiuin} and \eqref{ultimomomento}
 are consequences of this last fact
 (see the comments after Example \ref{comments before}).

Assume that majorization holds in \eqref{hora de comer}, so then we
have $$\tr(I_d^2)=\tr(\frac{d}{l} I_l)=\tr(\mathbf
P_q(V))=\tr((VV^*)^2)=\tr((V^*V)^2)$$ Since $I_d\prec_w V^*V$ and
the function $f(x)=x^2$ is strictly convex, by Theorem \ref{props
submayo} we conclude that there exists a unitary $U\in \U(d)$ such
that
$$V^*V=U^*(  I_d)U=I_d.$$ On the
other hand, if there exists an u.i.n. $\|\cdot\|$ such that equality
holds in \eqref{desiuin} then using the left-hand side of
\eqref{tempranillo} we get
 \begin{equation}
 d\cdot \eta_\psi(l)=\|\mathbf P_q(V)\|\geq \tr(\mathbf P_q(V))\cdot \eta_\psi(l)\geq  d\cdot \eta_\psi(l)
 \end{equation}which implies that $\tr((V^*V)^2)=\tr(\mathbf
 P_q(V))=d$. As before, we conclude that $V^*V=I_d$.
 Finally, it is clear that in case $\{V_i\}_i$ is a projective rank-$l$ $(m,l,d)$-protocol
 then $\mathbf P_q(V)=\frac{d}{l}I_l$. The last part of the theorem follows from this fact.
 \end{proof}

\begin{teor}[Optimality of u.w.p. protocols]\label{infi}
Let $\{V_i\}_{i=1}^m$ be an $(m,l,d)$-reconstruction system on the Hilbert
space $\mathcal H$ such that $\tr((V_i^*V_i)^{1/2})\geq (\frac{d\cdot
l}{m})^{1/2}$ for $1\leq i\leq m$. Then,
\begin{equation}\label{horas} \frac{d}{l}\,I_l\prec_w \mathbf P_q(V)
\end{equation}Hence,
for every u.i.n. $\|\cdot\|$ on $M_l(\CC)$ with associated symmetric gauge
function $\psi$ we have
\begin{equation}\label{desiuin2} d\cdot \eta_\psi(l)\leq \|\mathbf
P_q(V)\|
\end{equation}
and for every increasing convex function $f:\RR_{\geq 0}\rightarrow
\RR$ with $f(0)=0$ we have \begin{equation}\label{ultimomomento2}
l\cdot f(\frac{d}{l})\leq \tr(f(\mathbf P_q(V))).
\end{equation}

Moreover, majorization holds in \eqref{horas} or there exists u.i.n.
$\|\cdot\|$ such that equality holds in \eqref{desiuin2} or there
exists an increasing strictly convex function $f:\RR_{\geq
0}\rightarrow \RR$ with $f(0)=0$ such that equality holds in
\eqref{ultimomomento2}  if and only if $\{V_i\}_i$ is a u.w.p.
rank-$l$ $(m,l,d)$-protocol.
\end{teor}

\begin{proof}
Let $\{V_i\}_i$ be an $(m,l,d)$-reconstruction system such that, for
$1\leq i\leq m$
\begin{equation*}\label{empe}
\tr((V_i^*V_i)^{1/2})=\tr((V_iV_i^*)^{1/2})=(\frac{d\cdot l}{m})^{1/2} \ \
\Rightarrow (\frac{d}{m\,l})^{1/2}\, I_l\prec_ w (V_iV_i^*)^{1/2}
\end{equation*}  and thus
 $\tr(V_i^*V_i)=\tr(V_iV_i^*)\geq \tr( \frac{d}{m\,l}
 I_l)=\frac{d}{m}$. Hence, $\sum_{i=1}^m \tr(V_i^*V_i)\geq d$   and \eqref{horas} and \eqref{desiuin2} are
consequences of Theorem \ref{propopot}. If majorization holds in
\eqref{horas} or there exists u.i.n. $\|\cdot\|$ such that equality
holds in \eqref{desiuin2} or there
exists an increasing strictly convex function $f:\RR_{\geq
0}\rightarrow \RR$ with $f(0)=0$ such that equality holds in
\eqref{ultimomomento2} then again by Theorem \ref{propopot}, we
conclude that $\{V_i\}$ is an $(m,l,d)$-protocol. Thus, $\mathbf
P_q(V)=\sum_{i=1}^m V_iV_i^*$ with $\tr(\mathbf P_q(V))=d$.
Therefore, $\tr(V_iV_i^*)=\frac{d}{m}$ and since
$(\frac{d}{m\,l})^{1/2}\, I_l\prec_ w (V_iV_i^*)^{1/2}$ (recall that
$f(x)=x^2$ is an strictly convex function) we conclude as before
that $V_i^*V_i=\frac{d}{m\,l} P_i$ for some rank-$l$ orthogonal
projection $P_i$ for $1\leq i\leq m$.
\end{proof}

There are other issues regarding this potential, such as the structure of local minimizers where we consider the relativization of the product topology in $\prod_{i=1}^m \mathcal L(\mathcal H,\mathcal K)$, to the sets of reconstruction systems considered in the previous theorems. We shall consider these and related problems elsewhere.
\section{Optimal protocols for erasures and strict compatible u.i.n.}\label{sec 4}

Following \cite{Bod}(see also \cite{BodPau}, \cite{HolPau}) we begin
by modeling the situation in which in an encoding-decoding scheme
based on an $(m,l,d)$-protocol some fixed number of packets $(V_ix)$
are lost, corrupted or just delayed for such a long time the we
decide to reconstruct $x$ without these packets.

In  order to model the previous situation we consider a signal as a
vector in the $d$-dimensional (real or complex) vector space
$\mathcal H$, which is transmitted in the form of $m$ packets of $l$
linear coefficients. Hence, each packet is a vector the
$l$-dimensional (real or complex) Hilbert space $\mathcal K$.
 We shall assume that
$d< ml$ to allow for redundancy of the information sent through the
channel and thus for the possibility of a reasonable reconstruction
even when some packets are lost in the transmission. On the other
hand, we shall also assume that $l<d$ i.e. the dimension
(complexity) of the data is strictly bigger that the dimension of
the noiseless sub-channel (sub-system) which constitute the packets
(otherwise there are trivial optimal schemes).

Given $\mathbb K\subseteq \mathbb J:=\{1,\ldots,m\}$ a subset of
size $|\mathbb K|=p$ we consider the associated \emph{packet-lost
operator} $E_\mathbb K$ on $\oplus_{j=1}^m\mathcal K$ given by
$E_\mathbb K(\oplus_{j=1}^m y_j)=\oplus_{j=1}^m (1-\chi_\mathbb
K(i))\, y_i$, where $\chi_\mathbb K:\mathbb J\rightarrow \{0,1\}$
denotes the characteristic function of the set $\mathbb
K\subset\mathbb J$. We denote $D_\mathbb K:=I-E_\mathbb K$. In order
to simplify the notation we write $E_j$ (respectively $D_j$) in case
$\mathbb K=\{j\}$.

In our present situation, we shall consider a ``blind
reconstruction'' strategy for $(m,l,d)$-protocols for $\mathcal H$.
In case some packets are lost, i.e. assuming that the encoded
information $Vx\in \oplus_{i=1}^m\mathcal K$ (for some $x\in
\mathcal H$) is altered according to the packet-lost operator
$E_\mathbb K$, our reconstructed vector will be $V^*E_\mathbb
KV(x)$, where $V$ denotes the analysis operator of the
$(m,l,d)$-protocol $\{V_i\}_{i=1}^m$.

As a measure of performance of an $(m,l,d)$-protocol in this setting
we introduce the worst-case reconstruction error when $p$ packets are lost with
respect to an arbitrary compatible unitarily invariant
norm:\begin{equation*}\label{defi de error}
e_p^\psi(V):=\max\{\|V^*V-V^*E_\mathbb K V\|: \ \mathbb K\subseteq
\mathbb J, \ |\mathbb K|=p \}
\end{equation*}where $\|\cdot \|$ is a compatible u.i.n.
with associated symmetric gauge function $\psi$  (see Definition
\ref{defi de uin comp}) and $V$ denotes the analysis operator of the
$(m,l,d)$-protocol $\{V_i\}_{i=1}^m$.  Since the set $\mathcal
V(m,l,d)$ of all $(m,l,d)$-protocols is compact then the value
\begin{equation*} e_1^\psi(m,l,d)=\inf \{e_1^\psi(V):\ \{V_i\}_i\in
\mathcal V(m,l,d)\}
\end{equation*} is attained and we define the set of \emph{1-loss
optimal protocols for $\|\cdot\|$} to be the nonempty compact set
$\mathcal V_1^\psi(m,l,d)$ where this infimum is attained, i.e.
\begin{equation*}
\mathcal V_1^\psi(m,l,d):=\{\{V_i\}_i\in \mathcal V(m,l,d):
e_1^\psi(V)=e_1^\psi(m,l,d)\}
\end{equation*} Proceeding inductively, we now set for $1\leq p\leq
m$ \begin{equation*} e_p^\psi(m,l,d)=\inf\{e_p^\psi(V):\ \{V_i\}_i\in
\mathcal V_{p-1}^\psi(m,l,d)\}
\end{equation*}and define the \emph{optimal $p$-protocols for
$\|\cdot\|$} to be the non-empty compact subset of $\mathcal
V_{p-1}^\psi(m,l,d)$ where this infimum is attained.

\subsection{$\mathbf{e_1^\psi(\cdot)}$ optimality for one package lost}\label{subsec 4.1}

\begin{lema}
Let $\|\cdot\|$ be a compatible u.i.n. with associated symmetric gauge
function $\psi$. Let $\{V_i\}_{i=1}^m$ be an $(m,l,d)$-protocol on the Hilbert space $\mathcal
H$. Then,
\begin{equation}\label{desienun1}\max_{1\leq j\leq m}
\|V_jV_j^*\|\geq \frac{d\cdot \eta_\psi(l)}{m},\end{equation} where
$\eta_\psi(l)=\frac{\psi(\mathbf e_l)}{l}$. Moreover, if $\|\cdot
\|$ is strict then equality holds in \eqref{desienun1} if and only
if $\{V_i\}_{i=1}^m$ is a u.w.p. rank-$l$ $(m,l,d)$-protocol.
\end{lema}

\begin{proof}
Following \cite{Bod} we consider \begin{equation}\label{desi1}
\max_j \|\, V_j\, V_j^* \, \|\geq \frac{1}{m}\, \sum_{j=1}^m \|\,
V_j\, V_j^* \, \|.\end{equation} Recall that in this case
$\frac{\tr(V_j\,V_j^*)}{l}\ \mathbf e_l\prec
\lambda(V_j\,V_j^*)$
and hence
\begin{equation}\label{desi2}
\|V_j\,V_j^*\|\geq \frac{\tr(V_j\,V_j^*)}{l}\ \psi(\mathbf
e_l)=\tr(V_jV_j^*)\,\eta_\psi(l).\end{equation}Using the fact that
$\sum_{i=1}^m\tr(V_j\,V_j^*)=d$, \eqref{desienun1} now follows from
\eqref{desi1} and \eqref{desi2}.

Assume further that $\|\cdot\|$ is strict and the equality holds in
\eqref{desienun1}. Then, equality also hold in \eqref{desi1} and
\eqref{desi2}, too. Since $\|\cdot\|$ is strict  we conclude that
$\lambda(V_jV_j^*)=\frac{\tr(V_jV_j^*)}{l}\,\mathbf e_l$ and hence
$V_j^*V_j$ is a multiple (independent of $j$) of a rank-$l$ projection. The     lemma
easily follows from these facts.
\end{proof}

\begin{teor}
Let $\|\cdot\|$ be a compatible u.i.n.g with associated symmetric gauge
function $\psi$. Let $\{V_i\}_{i=1}^m$ be the coordinate
operators of an $(m,l,d)$-protocol on the Hilbert space $\mathcal
H$. Then,
\begin{equation}\label{desiteo1}
e_1^\psi(V)\geq \frac{d\cdot \eta_\psi(l)}{m}.
\end{equation}
Moreover, if $\|\cdot\|$ is strict then equality holds in
\eqref{desiteo1} if and only if $\{V_i\}_{i=1}^m$ is a u.w.p. rank-$l$ $(m,l,d)$-protocol.
\end{teor}

\begin{proof}
 For fixed $1\leq j\leq m$ note
that $V^*V-V^*E_jV=V^*D_j\,V$ and
$$\|V^*D_jV\|=\|D_jVV^*D_j\|=\|V_jV_j^*\|=\|V_j^*V_j\|.$$ Therefore,
the quantity to be minimized is $e_1^\psi(V)=\max_j\|V_j^*V_j\|$.
The result now follows from the previous lemma.
\end{proof}

The previous theorem completely characterizes the structure of the 1-loss optimal $(m,l,d)$-protocols in case $\|\cdot\|$ is an strict compatible u.i.n. Since the operator norm is a compatible strict u.i.n. we derive in particular \cite[Theorem 13]{Bod} (note that for the operator norm $\|\cdot\|_{\infty}$ we have $\eta_{\infty}(l)=\frac{1}{l}$). In section \ref{sec 5} we shall be concerned with the existence of protocols with prescribed properties (such as u.w.p. rank-$l$ $(m,l,d)$-protocols).

\subsection{The case of two lost packages}\label{subsec 4.2}

Consider the quantity defined in \cite{Bod} $$
c_{\,m,\,l,d}:=\sqrt{\frac{d}{(m-1) \ m\, l  }\ (1-\frac{d}{m l})}.
$$

In what follows we consider the class
\begin{equation*}\label{defi de c}
\mathcal  C(m,l,d)=\{ \{V_i\}_i: \ \text{u.w.p. rank-}l\  (m,l,d) \,\text{
protocol, } \max_{1\leq i\neq j\leq m}\tr(|V_iV_j ^*|)\geq l\cdot
c_{\,m,\,l,d}\}
\end{equation*}

\begin{teor}[${e_2^\psi}$ optimality in ${\mathcal C(m,l,d)}$]\label{elteo}
Let $\|\cdot\|$ be a compatible u.i.n. with associated symmetric gauge
 function $\psi$. Then, if $\{V_i\}_i\in \mathcal C(m,l,d)$
we have that
\begin{equation}\label{ecteo22}
e_2^\psi(V)\geq \psi((\frac{d}{ml}+c_{\,m,\,l,d})\ \mathbf e_l,(\frac{d}{ml} - c_{\,m,\,l,d})\ \mathbf e_l)
\end{equation} If $\{V_i\}_i$ is a u.w.p rank-$l$ $(m,l,d)$ protocol such
that for $i\neq j$ $V_iV_j^*=c_{\,m,\,l,d}\, Q_{i,j}$ for unitary
operators on $\mathcal K$, then $\{V_i\}_i\in \mathcal C(m,l,d)$ and
it attains the bound for $e_2^\psi$ in \eqref{ecteo22}.
\end{teor}
\begin{proof}
In order to compute the worst case reconstruction error for two lost
packages we note that if $\|\cdot\|$ is a compatible u.i.n. then (see the comments after Definition \ref{defi
de uin comp} in the Preliminaries)
\begin{eqnarray*} \| V^*(D_i+D_j)V \|&=&\| (D_i+D_j) V\,V^*(D_i+D_j) \|
 =\left\| \ \begin{pmatrix}
\frac{d}{m l} \,I & V_iV_j^*
\\ V_jV_i^* & \frac{d}{m l} \,I
\end{pmatrix} \ \right\|\\
&=&\psi((\frac{d}{m\, l} \, \mathbf e_l + s(V_iV_j^*),\frac{d}{m\,
l} \, \mathbf e_l - s(V_iV_j^*)))
\end{eqnarray*} where the last equality above follows from \cite[Theorem
7.3.7]{HoJo} and $s(A)=\lambda(|A|)\in \RR^l$ is the vector of
singular values of $A\in M_l(\CC)$. Notice that for $i\neq j$
\begin{equation}\label{igualtr}
\tr((\frac{d}{m\, l} \, \mathbf e_l + s(V_iV_j^*),\frac{d}{m\,
l} \, \mathbf e_l - s(V_iV_j^*)))=2\frac{d}{m} \ ,
\end{equation}
and since $\{V_i\}_i\in \mathcal C(m,l,d)$ then for some fixed
$i_0\neq j_0$ we should have
\begin{equation}\label{cotatr}\tr((\frac{d}{m\, l} \, \mathbf e_l +
s(V_{i_0}V_{j_0}^*))=\frac{d}{m}+\tr(|V_{i_0}V_{j_0}^*|)\geq
\frac{d}{m }+l\cdot c_{\,m,\,l,d}.\end{equation} Now,
\eqref{igualtr}, \eqref{cotatr} and Lemma \ref{lemmayo} imply that
in this case \begin{equation*} ((\frac{d}{m l}+c_{\,m,\,l,d})\
\mathbf e_l, (\frac{d}{m l}-c_{\,m,\,l,d})\ \mathbf e_l)\prec
(\frac{d}{m\, l} \, \mathbf e_l + s(V_{i_0}V_{j_0}^*),\frac{d}{m\,
l} \, \mathbf e_l - s(V_{i_0}V_{j_0}^*))\, .  \end{equation*}
Therefore,
\begin{equation*}\label{serepite}e_2^\psi(V)\geq \|V^*(D_{i_0}+D_{j_0})V\|\geq
\psi((\frac{d}{ml}+c_{\,m,\,l,d})\ \mathbf e_l,(\frac{d}{ml} -
c_{\,m,\,l,d})\ \mathbf e_l)\, .\end{equation*} Finally, is clear
that in case that $\{V_i\}_i$ is such that for $i\neq j$,
$V_iV_j^*=c_{\,m,\,l,d}\, Q_{i,j}$ for unitary operators on
$\mathcal K$,
 then $ \{V_i\}_i \in \mathcal C(m,l,d)$  and it attains the bound of $e_2^\psi$ in
\eqref{ecteo22}.

\end{proof}

It would be interesting to characterize the structure of all u.w.p. rank-$l$
$(m,l,d)$-protocols that attain the lower bound in \eqref{ecteo22}
in the general context of compatible u.i.n. On the other hand, it is
not clear at this point whether the condition in the definition of
the class $\mathcal C(m,l,d)$ is not trivial, i.e. it holds for
every u.w.p. protocol (see also Lemma \ref{truch} and Corollary
\ref{coroa}).

\medskip

The following facts are known for $l=1$ (see \cite{HolPau}).

\medskip

\begin{lema}\label{truch}
Let $\|\cdot\|$ be a compatible u.i.n.
with associated symmetric gauge function $\psi$. Let $\{V_i\}_{i=1}^m$ be a u.w.p.
rank-$l$ $(m,l,d)$-protocol on the Hilbert space $\mathcal H$.
Then, for every $1\leq i\leq m$ we have
\begin{equation}\label{iguanorcuad}
\sum_{j=1,\ j\neq i}^m \tr(|V_jV_i^*|^2)=\frac{d}{m}\,(1- \frac{d}{m l} ),
\end{equation}
\begin{equation}\label{desicota}
 \sum_{j=1,\ j\neq i}^m \tr(|V_iV_j^*|)\geq \sqrt{\frac{d}{ml}\,(1-
\frac{d}{m l} )}\ \cdot l
\end{equation}
and hence
\begin{equation}\label{ecanexa}
\max_{1\leq j\leq m,\ i\neq j} \ \tr(|V_iV_j^*|^2)\geq c^2_{m,l,d} \cdot l,
\end{equation}

\begin{equation}\label{otro truch}
\max_{1\leq j\leq m,\ i\neq j}\tr(| V_i V_j^*|)\geq
\max\left(\frac{c_{\,m,\,l,d}\cdot l}{\sqrt{m-1}}\, , \, c_{\,m,\,l,d}\cdot
\sqrt{l}\right).
\end{equation}

\end{lema}

\begin{proof}
Since $VV^*=(VV^*)^2$ then for fixed $1\leq i\leq m$
\begin{equation} \label{yoto}\frac{d}{m l}\,I_l=V_iV_i^*=\sum_{j=1}^m
|V_jV_i^*|^2=\sum_{j=1,\ j\neq i}^m  |V_jV_i^*|^2+ \frac{d^2}{m^2 l^2}\,I_l.
\end{equation}  Now
\eqref{iguanorcuad} follows by taking traces in \eqref{yoto}. Using
again \eqref{yoto} and the concavity of the square root function
\cite{Rot} we get
\begin{equation*}
\sum_{j=1,\ j\neq i}^m  \tr(|V_jV_i^*|)\geq \tr(\sqrt{ ( \frac{d}{m
l}-\frac{d^2}{m^2 l^2})I_l} )
\end{equation*} which is \eqref{desicota}. Now, from
\eqref{iguanorcuad} we get \eqref{ecanexa}.
Using \eqref{ecanexa} we get that, for fixed $1\leq i\leq m$  \begin{equation}\label{uf}
\max_{1\leq j\leq m,\ i\neq j} \tr(|V_iV_j^*|)\geq \max_{1\leq j\leq m,\ i\neq j} \tr(|V_iV_j^*|^2)^{1/2}\geq
\sqrt{c_{\,m,\,l,d}^2 \cdot l}.
\end{equation}Finally, from \eqref{desicota} and using \eqref{uf} we get
 \eqref{otro truch}.
\end{proof}

\begin{rema}\label{sincoment}
Under the hypothesis of Lemma \ref{truch}, note that \eqref{ecanexa}
implies that, for fixed $1\leq i\leq m$ then
\begin{eqnarray*}\max_{1\leq j\leq m, \, i\neq
j}\|\, |V_iV_j^*|^2 \| \geq \frac{1}{m-1}\sum_{1\leq j\leq m, \, i\neq
j} \| \,|V_iV_j^*|^2 \|  \geq \\ \nonumber c_{\,m,\,l,d}^2(\psi):=\frac{1}{m-1}\sum_{1\leq j\leq m, \, i\neq
j}\tr( |V_iV_j^*|^2) \, \eta_\psi(l)
=\frac{d\cdot
\eta_\psi(l)}{m\ (m-1)}\ (1-\frac{d}{m l}).
\end{eqnarray*}
If we assume further that $\|\cdot\|$ is strict and that for fixed $1\leq i\leq m$ $$ \max_{1\leq j\leq m, \, i\neq
j}\|\, |V_iV_j^*|^2 \| = c_{\,m,\,l,d}^2(\psi)$$
 then for every $j\neq i$, $|V_iV_j^*|$
has only one eigenvalue, namely $c_{\,m,\,l,d}$. Using the polar
decomposition for $V_iV_j^*$ we conclude that
$V_iV_j^*=c_{\,m,\,l,d} \,Q_{i,j}$ for some unitary operator
$Q_{i,j}$ in $\mathcal K$. In particular, \begin{equation}\label{desilem3}  \max_{1\leq i\neq j\leq
m}\|\, |V_iV_j^*|^2 \| \geq c_{\,m,\,l,d}^2(\psi)\end{equation} and equality
holds if and only if, for $1\leq i\neq j\leq m$ then
$V_iV_j^*=c_{\,m,\,l,d}\,Q_{ij}$ for unitary operators $Q_{ij}$ in
$\mathcal K$. These remarks generalize to this context \cite[Lemma
14]{Bod} for the spectral norm (notice that in this case
$\eta_\infty(l)=\frac{1}{l}$); in particular, \eqref{desilem3} is an
extension of a result of Welch \cite{wel}.
\end{rema}
\smallskip

Given a compatible strict u.i.n. $\|\cdot\|$ we say that it is $k$-strongly
strict if for every $A,\,B\in M_k(\CC)^{sa}$ such that $A\prec B$
and $\|A\|=\|B\|$ then $A=U^*BU$ for some $U\in \U(k)$. For example, the $p$-norms are $k$-strongly strict for $k\geq 1$
(see Example \ref{comments before}). On the other hand, it is easy
to see that the operator norm is 2-strongly strict.

\begin{teor}\label{coroa}
Let $\|\cdot\|$ be a compatible u.i.n. with associated symmetric gauge
 function $\psi$.
\begin{enumerate}
\item If $\{V_i\}_i$ is a u.w.p. $(m,1,d)$-protocol (i.e. a uniform tight
frame of $m$ vectors) then $\{V_i\}_i\in \mathcal C(m,1,d)$ and
\begin{equation}\label{operadora}
e_2^\psi(V)\geq \psi((\frac{d}{m}+c_{\,m,1,d},\frac{d}{m} - c_{\,m,\,1,d}))
\end{equation}
If we further have that
$V_iV_j^*=c_{\,m,\,1,d}\ q_{ij}$ for $q_{ij}\in \CC$ with
$|q_{ij}|=1$, for every $i\neq j$ then equality holds in \eqref{operadora}. Moreover, the converse is true for 2-strongly strict compatible u.i.n. 
\item If $\{V_i\}_i$ is a u.w.p. rank-$l$ $(2,l,d)$-protocol
then $\{V_i\}_i\in \mathcal C(2,l,d)$ and \begin{equation*}
e_2^\psi(V)\geq \psi((\frac{d}{2l}+c_{\,2,l,d})\,\mathbf
e_l,(\frac{d}{2l} - c_{\,2,\,l,d})\,\mathbf e_l)
\end{equation*} If $\{V_i\}_i$ is a u.w.p-$(2,l,d)$ protocol such
that for $i\neq j$, $V_iV_j^*=c_{\,2,\,l,d}\,Q_{ij}$ for unitary
operators $Q_{ij}$ in $\mathcal K$, it attains the bound for
$e_2^\psi$ above.

\end{enumerate}
\end{teor}
\begin{proof}
By setting respectively $l=1$, respectively $m=2$,  in \eqref{otro
truch} we see that in these cases $\mathcal C(m,l,d)$ coincides with
the class of all u.w.p. rank-$l$ $(m,l,d)$-protocols (i.e. the
condition in the definition of $\mathcal C(m,l,d)$ becomes trivial
in these cases) so the first part of item (i) and item (ii) follow
from Theorem \ref{elteo}.

In order to prove the second assertion in item (i) assume that $\|\cdot\|$ is a 2-strongly strict compatible u.i.n. Note that if
$\alpha,\,\beta\in \RR^2$ are such that $\tr(\alpha)=\tr(\beta)$
then these vectors are comparable with respect to majorization; indeed $\alpha\prec \beta$ if and only if $\max\{\alpha_1,\alpha_2\}\leq \max\{\beta_1,\beta_2\}$. Assume now that $\|\cdot\|$ is a 2-strongly strict norm and that $\{V_i\}_i$ is an u.w.p. $(m,1,d)$-protocol in which the lower bound in \eqref{operadora} is attained. Hence, by inspection  of the
proof of Theorem \ref{elteo} (note that $V_iV_j^*\in \CC$ for $l=1$) we see that if $i\neq j$ then
$$ \psi(\frac{d}{m} + |V_iV_j^*|,
  \frac{d}{m} - |V_iV_j^*|)\leq \psi(\frac{d}{m} + c_{m,1,d}, \frac{d}{m} - c_{m,1,d})$$
  which implies that
  \begin{equation}\label{implicado} \frac{d}{m} + |V_iV_j^*|
  \leq \frac{d}{m} + c_{m,1,d} \ \ \Rightarrow |V_iV_j^*|
  \leq c_{m,1,d} \ , \ i\neq j.\end{equation}
  Since $$ \tr(VV^*)=\tr((VV^*)^2)=\sum_{i\neq j}|V_iV_j^*|^2+\frac{d^2}{m}=\sum_{i\neq j}c_{m,1,d}^2+\frac{d^2}{m}$$
  we conclude that equality holds
  in the right hand side of \eqref{implicado} and the theorem follows from this last fact.
\end{proof}

\begin{rema}\label{esppol}The first item in Theorem \ref{coroa}
generalizes the results in \cite{BodPau} and \cite{HolPau} about the
optimality of 2-uniform frames to the context of strongly strict compatible unitarily invariant
norms.
\end{rema}

\section{Existence of optimal protocols for one package lost and the $q$-fundamental inequality}\label{sec 5}

In \cite{Bod}, \cite{BodPau}, \cite{HolPau}, \cite{Pau} and
\cite{Tropp}, several examples of 2-loss optimal protocols, i.e
u.w.p. rank-$l$ $(m,l,d)$-protocols $\{V_i\}$ for which
$V_iV_j^*=c_{m,l,d}\,Q_{ij}$ with $Q_{ij}\in \U(l)$, are constructed
based on different techniques. Still, the problem of finding
necessary and sufficient conditions for the existence of 1-loss
optimal protocols, i.e. u.w.p. rank-$l$ $(m,l,d)$-protocols, has
been considered open (see the discussion in \cite{Pau} about this
topic).

In the case $l=1$ (i.e. the classical case of frames), the existence
of tight normalized frames with given norms of the frame vectors
(and hence of 1-loss optimal protocols) is characterized completely
by the so-called \emph{fundamental frame inequality} discovered in
\cite{casazza2}. Moreover it is now known (\cite{JDM},
\cite{casazza3}, \cite{Lar}, \cite{MR}) that the fundamental frame
inequality is a particular case of a majorization relation (via the
Schur-Horn theorem) that constitutes a necessary and sufficient
condition for the existence of a frame with prescribed norms of the
frame vectors and frame operator.

In what follows we exhibit necessary and  sufficient (spectral)
conditions for the existence of $(m,l,d)$-protocols $\{V_i\}_i$ with
prescribed eigenvalue vectors $\lambda(V_i^*V_i)\in \RR_{\geq 0}^d$
for $1\leq i\leq m$. As in the classical case $l=1$ there exists a
relation between these conditions and an extended notion of (block)
majorization as introduced in \cite{Mas}.

\begin{teor}\label{con suf y nec}
Let $\lambda_i=\lambda_i^\downarrow\in \RR_{\geq 0}^l$ for $1\leq
i\leq m$. Then, the following statements are equivalent:
\begin{enumerate}
\item\label{item 1} There exists an $(m,l,d)$-protocol $\{V_i\}_{i=1}^m$ such that  $\lambda(V_i^*V_i)=(\lambda_i, 0_{d-l})$, for $1\leq i\leq m$.
\item\label{item 2} There exist
$\{A_i\}_{i=1}^m\subset M_d(\FF)^+$ such that
$$\lambda(A_i)=(\lambda_i, 0_{d-l})\ \text{ for }\ 1\leq i\leq m \ \ \text{ and } \ \ \sum_{i=1}^mA_i=I_d.$$
\item\label{item 3} The $(m+1)$-tuple $$((\lambda_1,0_{d-l}),
\ldots,(\lambda_m,0_{d-l}),\mathbf e)\in (\RR^d)^{(m+1)}$$
 satisfy
Klyachko's compatibility inequalities plus $\sum_{i=1}^m\tr(\lambda_i)=d$.
\item\label{item 4} There exists an orthogonal projection $P\in M_m(M_l(\FF))$
 with $\tr(P)=d$ and such that, if $P=(P_{ij})_{i,j=1}^m$ with $P_{ij}\in M_l(\FF)$
 for $1\leq i,j\leq m$, then
 $$\lambda(P_{ii})=\lambda_i \ \text{ for } \ 1\leq i\leq m.$$
\end{enumerate}
\end{teor}

\begin{proof}

Clearly, (\ref{item 1}) implies (\ref{item 2}) by considering $A_i=V_i^*V_i$ for $1\leq i\leq m$. Assume
then item (\ref{item 2}). In this case note that rank$(A_i)\leq l$ and hence
there exist linear operators $V_i:\mathcal H\rightarrow \mathcal K$
such that $V_i^*V_i=A_i$ for $1\leq i\leq m$. It is clear that $\{V_i\}_{i=1}^m$ is an $(m,l,d)$-protocol as in (\ref{item 1}). Therefore, (\ref{item 1}) and (\ref{item 2}) are equivalent.

The equivalence of items (\ref{item 2}) and (\ref{item 3}) is Theorem \ref{teoK}.

Assume again (\ref{item 1}) holds and let $V:\mathcal H\rightarrow \oplus_{i=1}^m\mathcal K$ be the analysis operator of the protocol $\{V_i\}_i$.
 Since $
V^*V=\sum_{i=1}^mV_i^*V_i=1_d$ then we get that the block matrix $
VV^*=(V_iV_j^*)_{i,j=1}^m\in M_m(M_l(\FF))$ (i.e. the Grammian of
$\{V_i\}_i$) is an orthogonal projection; moreover, note that
$\tr(VV^*)=\tr(V^*V)=d$ and that the diagonal blocks of the Grammian
satisfy $(\lambda(V_iV_i^*),0_{d-l})=\lambda(V_i^*V_i)=(\lambda_i,
0_{d-l})$, for $1\leq i\leq m$ (see the comments after Definition
\ref{defi de uin comp}). Conversely, assume that item (\ref{item 4})
holds and let $V:\mathcal H\rightarrow \oplus_{i=1}^m\mathcal K$ be
an isometry such that $VV^*=P$ (such an isometry exists since
rank$(P)=d$ by assumption). Let $V_i:\mathcal H\rightarrow \mathcal
K$ for $1\leq i\leq m$ be such that $Vx=\oplus_{i=1}^m V_ix$ and
note that then $P=VV^*=(V_iV_j^*)_{ij}$ and that
$I_d=V^*V=\sum_{i=1}^mV_i^*V_i$ that is, $\{V_i\}_{i}$ is an
$(m,l,d)$-protocol as in (\ref{item 1}). Thus, items (\ref{item 1})
and (\ref{item 4}) are equivalent.
\end{proof}

\begin{rema}\label{comentarios LiPoon} Using the characterization in
item (\ref{item 4}) in Theorem \ref{con suf y nec} and the reduction
described in \cite{LiPoon} (which is relevant from an algorithmic
point of view) it is possible to show that Klyachko's compatibility
inequalities in \eqref{item 3} in Theorem \ref{con suf y nec}
 can be reduced to a system of inequalities that, in case $l=1$
are simply the conditions given in the majorization relation
diag$(\|P_{11}\|^2,\ldots,\|P_{mm}\|^2)\prec I_d\,\oplus \, \mathbf
0_{m\cdot l-\,d}$, where diag$(x)\in M_n(\CC)$ is the diagonal
matrix with main diagonal $x\in \CC^n$.

Actually, the inequalities in \eqref{item 3} in Theorem \ref{con suf
y nec} can be regarded as determining an extended notion of
majorization as defined in \cite{Mas}. Indeed, with the terminology
of \cite[Definition 4.4]{Mas}, the conditions given in Theorem \ref{con
suf y nec} are also equivalent to the $\mathbf t$-extended
majorization relation $\oplus_{i=1}^m
\text{diag}(\lambda_i)\prec_{\mathbf t} I_d \,\oplus \, \mathbf
0_{m\cdot l-\,d}\in M_{m\cdot l}(\CC)$, where $t=(\mathbf
e_l,1)_{i=1}^m$.
\end{rema}

\begin{coro}[$q$-fundamental projective $(m,l,d)$-protocol
inequalities]\label{q fund ineq} Let $t(i)\in \{1,\ldots,l\}$ and
$w_i\in \RR_{\geq 0}$ for $1\leq i\leq m$. Then, there exists a
projective $(m,l,d)$-protocol $\{V_i\}_{i=1}^m$ for the Hilbert
space $\mathcal H$ such that $V_i^*V_i=w_i\,P_i$ for orthogonal
projections $P_i$ with $\tr(P_i)=t(i)$ for $1\leq i\leq m$
 if and only if for every $1\leq r\leq d$ and every
$(J_0,\ldots,J_m)\in LR_r^{\,d}(m)$ we have that
\begin{equation*}\label{la ec posta}
r\leq \sum_{i=1}^m w_i\cdot |\,J_i\cap \{1,\ldots, t(j)\} \,|
\end{equation*}plus the condition $d=\sum_{i=1}^m w_i\cdot t(i)$.
\end{coro}

As an immediate consequence of the $q$-f.p.p.i. we conclude that
u.w.p. rank-$l$ $(m,l,d)$-protocols exist if and only if for every
$1\leq r\leq d$ and every $(J_0,\ldots,J_m)\in LR_r^{\,d}(m)$ it
holds that
 $$ r\leq \frac{d}{m\cdot l} \cdot  \sum_{i=1}^m |\, J_i\cap \{1,\ldots,l\}\, |.$$

\begin{exa} Next, we show explicitly how to construct a
 projection $P=(P_{ij})_{ij}\in M_m(M_l(\CC))$ such that
 $P_{ii}=\frac{d}{ml}\,I$ for $1\leq i\leq m$, when
 $d=k\cdot l$ for some $k\in \NN$. Thus, by Theorem \ref{con suf y
 nec}, we show the existence of u.w.p. rank-$l$ $(m,l,d)$-protocols
 in this case.  This construction is
 a particular case of that appearing in the proof of \cite[Prop. 4.12]{Mas}: consider first
$\xi\in \CC$ an $m$-th primitive root of unity and let $\tilde U\in
\mathcal M_m(\CC)$ be the matrix with $j$-th row given by
$$R_j(\tilde
U)=1/\sqrt{m}\,(1,\,\xi^j,\,\xi^{2j},\ldots,\xi^{(m-1)j})\,,\ \
1\leq j\leq m.$$ It is then straightforward to show that the rows of
$\tilde U$ form an orthonormal basis for $\CC^m$ and hence $\tilde
U\in \U(m)$ is a unitary matrix. Let $U\in \U(d\cdot m)$ be the
block matrix $U=(\tilde U_{ij}\cdot 1_d)_{i,j=1}^m$. Then, consider the matrix $A=\oplus_{i=1}^k I=(A_{ij})_{ij}\oplus \mathbf 0_{(m-k)l}\in M_m(M_l(\CC))$ and note that
$$ U^*A U=(P_{ij})_{ij} ,\ \ P_{ii}=\frac{1}{m}\sum_{i=1}^m A_{ii}=\frac{k}{m}\,I \, ,$$
where the last equality follows from the diagonal block structure of
$A$ and by construction of $U$. Now, recall that $k=\frac{d}{l}$ and
we are done.
\end{exa}

\medskip

\noindent Pedro G. Massey \\
Departamento de Matem\'atica - FCE - Universidad Nacional de La
Plata\\P.O. Box 172 - Argentina\\
and\\
Instituto Argentino de Matem\'atica - CONICET - Argentina\\ e-mail
address: massey@mate.unlp.edu.ar

\end{document}